\documentclass[a4paper,11pt]{article}
\usepackage{amsfonts,amssymb,latexsym,amsmath}
\usepackage{multicol}
\begin{document}
\Large

\newcounter{num}[section]
\setcounter{num}{0}
\renewcommand{\thenum}{\arabic{section}.\arabic{num}}

\newcommand{\Num}{\refstepcounter{num}%
\textbf{\arabic{section}.\arabic{num}}}

\newcommand{\Theorem}{\textbf{Theorem~}}
\newcommand{\Prop}{\textbf{Proposition~}}
\newcommand{\Proof}{\textbf{Proof}}
\newcommand{\Ex}{\textbf{Example~}}
\newcommand{\Remark}{\textbf{Remark~}}
\newcommand{\Lemma}{\textbf{Lemma~}}
\newcommand{\Def}{\textbf{Definition~}}
\newcommand{\Cor}{\textbf{Corollary~}}
\newcommand{\Notations}{\textbf{Notations}}

\newcommand{\al}{\alpha}
\newcommand{\te}{\theta}
\newcommand{\la}{\lambda}
\newcommand{\La}{\Lambda}
\newcommand{\laSa}{\lambda_{S,a}}
\newcommand{\laSas}{\lambda_{S',a'}}
\newcommand{\VaSa}{V_{S,a}(\la)}
\newcommand{\oOSa}{\overline{\Omega}_{S,a}(\la)}
\newcommand{\OSa}{\Omega_{S,a}(\la)}
\newcommand{\OSas}{\Omega_{S',a'}(\la)}
\newcommand{\IOSa}{{\Ic_{S,a}}}
\newcommand{\Ola}{\Omega(\la)}

\newcommand{\eps}{\varepsilon}

\newcommand{\Sc}{{\cal S}}
\newcommand{\Mc}{{\cal M}}
\newcommand{\Pc}{{\cal P}}
\newcommand{\Ic}{{\cal I}}
\newcommand{\Fc}{{\cal F}}
\newcommand{\Ac}{{\cal A}}

\newcommand{\gx}{{\mathfrak g}}
\newcommand{\ogx}{\overline{\gx}}
\newcommand{\lx}{{\mathfrak l}}
\newcommand{\px}{{\mathfrak p}}
\newcommand{\bx}{{\mathfrak b}}

\newcommand{\ax}{{\mathfrak a}}

\newcommand{\gl}{{{\mathfrak g}{\mathfrak l}}}
\newcommand{\ut}{{{\mathfrak u}{\mathfrak t}}}
\newcommand{\UT}{{{\mathrm U}{\mathrm T}}}

\newcommand{\Ad}{{\mathrm{Ad}}}
\newcommand{\ad}{{\mathrm{ad}}}
\newcommand{\ind}{{\mathrm{ind}}}
\newcommand{\Hc}{{\cal H}}

\newcommand{\spp}{{\mathrm{span}}}
\newcommand{\codim}{{\mathrm{codim}}}
\newcommand{\rank}{{\mathrm{rank}}}

\newcommand{\Cb}{{\Bbb C}}
\newcommand{\Sb}{{\Bbb S}}
\newcommand{\Xb}{{\Bbb X}}
\newcommand{\Fbq}{{\Bbb F}_q}
\newcommand{\oFbq}{\overline{{\Bbb F}}_q}

\renewcommand{\leq}{\leqslant}
\renewcommand{\geq}{\geqslant}

\title{Representations of Gelfand-Graev type for  the unitriangular
 group}
\author{A.N.Panov\thanks{The paper is supported by the RFBR grants 12-01-00070-a,
12-01-00137-a, 13-01-97000-Volga region-a}\\
{\it Samara State University}\\
 apanov@list.ru}
\date{}
 \maketitle

{\small {\bf Keywords:} representation of the unitriangular group,
the orbit method, Hecke algebra

 \begin{center}
{\bf Abstract}
 \end{center}

 We consider the analog of Gelfand-Graev
representations of the uniteriangular group. We obtain the
decomposition into the sum of irreducible representations,  prove
that these representations are multiplicity free,  calculate the
Hecke algebra.}

\section{\bf {Introduction and main definitions}}

 A representation of Gelfand-Graev is a representation of the group $\mathrm{GL}(n,\Fbq)$
 (more generally, of the finite Chevalley group)
 induced from  a nondegenerate  character of its maximal unipotent subgroup.
The main property of  these representations that they are
multiplicity free.   These representations   appeared first in the
papers \cite{G1,G2}. It was proved that the algebra of
$G$-endomorphisms (further referred as the Hecke algebra)
 of these representations is commutative; this is equivalent to the
property of being multiplicity free.  The basis of the Hecke algebra
was constructed in the paper \cite{G2}, for $\mathrm{GL}(n,\Fbq)$,
and, later, for finite Chevalley groups, in \cite{I,S}. The
Gelfand-Graev representations play an important role in the
representation theory; many  papers appear on properties of these
representations and their generalizations.

 In this
paper, we study  representations $V(\la)$ of the unitiangular group;
this representations are  analogs of Gelfand-Graev ones. We shall
prove that $V(\la)$ are also multiplicity free. We shall give a
complete description of all its  irreducible components $\VaSa$.
Following the orbit method (see \cite{K1,K2}), we associate the
coadjoint orbit  $\OSa$ with the irreducible component $\VaSa$. In
the paper, we find a canonical form $\laSa\in\OSa$ and generators of
the defining ideal of the orbit $\OSa$. In a sequel of the paper, we
find a basis of the Hecke algebra $\Hc(\la)$. The main results are
formulated in theorems
\ref{th1},~\ref{th2},~\ref{th3},~\ref{th4},~\ref{th5}.

Decomposition of the representation  $V(\la)$ into  a sum of
irreducible components  admits  the interpretation in terms of the
theory of basic  characters (basic representations)   developed by
C. Andr\'e (see, for instance,   \cite{A1, A2, A3, A4}). One can
consider the representation  $V(\la)$  as the induced representation
from the basic (precisely, regular irreducible) representation of
the unitriangular subgroup of size  $n-1$ ~\cite{A5}.
 The representation  $V(\la)$ decomposes into the  representations
 $V_S(\la)$ that in its turn decomposes into irreducible components $\VaSa$ (see (\ref{nnn})).
Notice that the representations  $V_S(\la)$ are basic or sums of the
basic representations.

 Firstly, we shall give
a remark on the orbit method. The orbit method appeared in 1962 in
the paper \cite{K1}. There was
 shown that there exists one to one correspondence between
irreducible representations of a connected nilpotent Lie group  and
its coadjoint orbits. Later, in \cite{Kzh}, it was proved that the
orbit method is also true for any unipotent group over a finite
field (see also \cite{A1,P1}). There are some requirements on the
characteristic of the field; since the matrix exponent is used in
the orbit method, the characteristic of the finite field have to be
 great enough. For the group $\UT(n,\Fbq)$, it is sufficient to put
$\mathrm{char}\,\Fbq \geq n-1$.

In this paper, the characteristic is arbitrary.  The irreducible
representations  $\VaSa$ are induced from characters of associative
polarizations  (see definition \ref{pol}). This enables to remove
matrix exponents from the process of construction of representations
\cite{A2,Sa}.

 Let $\Fbq$ be a finite field of $q$ elements.
The unitriangular group   $G=\UT(n,\Fbq)$ consists of all upper
triangular matrices of size  $n\times n$ with units on the diagonal
and entries from the field  $\Fbq$. We consider further that $n>2$.
Denote by $\gx$ a subspace of all upper triangular matrices with
zeros on the diagonal. It is obvious that  $G=E+\gx$, where $E$ is
the unity matrix. The subspace $\gx$ is an associative algebra with
respect to the matrix multiplication and, therefore, a Lie algebra.

We shall give the following definition: a root is a pair $(i,j)$ of
integers, where $1\leq i<j\leq n$. The partial operation of addition
is defined on the set of all roots $$R=\{(i,j):~~1\leq i<j\leq n\}$$
as follows: $(i,j)+(j,s) = (i,s)$.

 The set of all roots $R$ decomposes  into the subsets
$R=R_+\sqcup R_0\sqcup R_-$, where $$ R_+ = \left\{ (i,j):~~
i+j<n+1\right\},$$ $$ R_0 = \left\{ (i,j):~~ i+j = n+1\right\},$$ $$
R_- = \left\{ (i,j):~~ i+j>n+1\right\}.$$

Matrix unities  $\{E_\al:~~ \al\in R\}$ generate  a basis in the
algebra $\gx$. The algebra $\gx$ is a direct sum of the subspaces
$\gx = \gx_+\oplus \gx_0\oplus\gx_-$, where
$$\gx_\pm = \spp \{ E_\al:~~ \al\in R_\pm\},$$
$$\gx_0
= \spp \{ E_\al:~~ \al\in R_0\}.$$ The subspaces  $\gx_+,~ \gx_0,
~\gx_-$ are subalgebras of the associative algebra $\gx$. Then
$$G_\pm = E+\gx_\pm,\quad ~ G_0 = E + \gx_0$$ are subgroups in  $G$.
 Any
root  $\al\in R$  determines  the one-parameter subgroup
 $$\{x_\al(t)
= E + tE_\al,~ t\in\Fbq\}.$$
 Introduce the following notations: \\
 1)~ $k=[\frac{n-1}{2}]$. Then  $n=2k+1$, if $ n $ is odd, $n=2(k+1)$, if $ n $ is even;\\
2)~~$\eps=\left\{\begin{array}{l} 0, ~~\mbox{if}~ n~~\mbox{is~~odd},\\
 1, ~~\mbox{if}~ n~~\mbox{is~~even}.\end{array}\right.$

We say that a root from $R_+$ is a simple root, if it can't be
presented as a sum of two roots from $R_+$. The set of all simple
roots is a union of two subsets  $\Pi_0\cup\Pi$, where
$$\Pi_0 = \{(i,i+1):~~1\leq i\leq k\},$$
$$\Pi = \{(i,n-i):~~1\leq i\leq k\}.$$
Notice that  $$\Pi_0\cap\Pi = \left\{\begin{array}{l} \varnothing,
~\mbox{if} ~~ n=2(k+1),\\
(k,k+1), ~\mbox{if} ~~ n=2k+1.\end{array}\right.$$
\Def\Num\label{chi}. A character of an associative algebra  $\ax$ is
a linear form  on it that is zero on  $\ax^2$.

Notice that any character of an associative algebra is  a character
of it as a Lie algebra.  Any character of the associative algebra
$\gx_+$  is uniquely determined  by its values on
$\{E_\al:~~ \al\in \Pi_0\cup \Pi\}$.\\
\Def\Num\label{nondeg}.  A character $\la:\gx_+\to\Fbq$
 is nondegenerate, if  $\la(E_\gamma)\ne 0$ for any $\gamma\in
\Pi\setminus \Pi_0$.

Fix some nontrivial complex character of the additive group of the
field $\Fbq$ (i.e. homomorphism $\Fbq\to\Cb^*$). We shall denote
this character by  $e^x$, where ~$x\in\Fbq$.

Any character $\la$ of the associative algebra  $\gx_+$ determines a
complex character (i.e. one-dimensional representation)  $\xi_\la$
of the group $G_+$ by the formula
\begin{equation} \label{xi}\xi_\la(1+x) =
e^{\lambda(x)},\quad x\in\gx_+.\end{equation} \Def\Num. The
representation  $V(\lambda) = \ind(\xi_\la,G_+,G)$, where $\lambda$
is a nondegenerate character of  $\gx_+$,
is called a representation of Gelfand-Graev type.\\
 \Prop\Num\label{VV}. \emph{If
$\la,~\la'$ are nondegenerate characters of  $\gx_+$ that coincide
on  $\{E_\al:~~ \al\in \Pi\}$, then  $V(\la)\cong
V(\la')$}. \\
\Proof. Let  $v_0\in V(\la)$ be the generating vector  of the
induced representation, $g_+v_0 = \xi_\la(g_+)v_0$. For any root
$\al\in \Pi_0\setminus \Pi$, there exists  a unique root
$\beta(\al)\in R_0$ such that the sum  $\gamma(\al) =
\al+\beta(\al)$ is determined
 and belongs to   $\Pi\setminus\Pi_0$. Indeed, if
$\al=(i,i+1)$, then $\beta(\al)=(i+1,n-i)$ and
$\gamma(\al)=(i,n-i)$. The subspace, spanned by $E_{\beta(\al)}$,
where  $\al\in \Pi_0\setminus \Pi$, is an associative subalgebra
with zero multiplication. Therefore, for any two roots $\al,\al'\in
\Pi_0\setminus \Pi$, the elements of the corresponding one-parameter
subgroups $x_{\beta(\al)}(t)$ and $x_{\beta(\al')}(t')$ commutes.
Consider the element  \begin{equation}\label{go} g_0 = \prod
x_{\beta(\al)}(t_\al)\in G,\end{equation} where
   $\al$ is running through $\Pi_0\setminus \Pi$, and
$t_\al\in \Fbq$ is  a solution of the equation
$$\la(E_\al) + \la(E_{\gamma(\al)})t_\al = \la'(E_\al).$$ Using the equalities
$$x_\al(s)x_{\beta(\al)}(t) =
x_{\beta(\al)}(t)x_\al(s)x_{\gamma(\al)}(st),$$ we obtain $g_+g_0v_0
= \xi_{\la'}(g_+)g_0v_0$. $\Box$

\section{\bf {Associative polarizations}}

Let  $\ax$ be an arbitrary  nilpotent associative algebra over an
arbitrary field  $K$. Adjoin the unity element  $E$ to the algebra
$\ax$. Then $G=E+\ax$ is an unipotent group.  The algebra  $\ax$ is
a Lie algebra with respect to the commutator  $[x,y]=xy-yx$. Let
$\la\in\ax^*$. Recall that  a polarization for $\la$ is a Lie
subalgebra $\px$ of $\ax$ that is a maximal isotropic subspace with
respect to the skew symmetric bilinear form
$B_\la(x,y) = \la([x,y])$. \\
\Def\Num\label{pol}. An associative polarization  of  $\la\in\ax^*$
is a
polarization $\px$ that obeys the following conditions\\
i)~ $\px$ is an associative subalgebra of  $\ax$,\\
ii) ~ $\la(\px^2)=0$.

Is it true that any  $\la\in\ax^*$ has  an associative polarization?
In general, the answer is negative.\\
\Ex\Num. The associative algebra
$\ax=\left\{\left(\begin{array}{lll}
0&a&b\\0&0&a\\0&0&0\end{array}\right)\right\}$ is commutative. Then
it is  a commutative Lie algebra. Any linear form $\la$ on $\ax$ has
a unique polarization, which coincides with $\ax$ . If
$\la(E_{13})\ne 0$, then  $\ax$ is not an associative polarization.

Suppose that $\ax$ is an associative nilpotent algebra over $\Fbq$.
 If $\px$ is an associative polarization for
$\la\in\ax^*$, then  formula (\ref{xi}) defines a one-dimensional
complex representation  $\xi_\la$ of the group  $P=E+\px$. Denote by
 $M(\la)$ the induced representation  $\ind(\xi_\la,P,G)$.
 \\
\Prop\Num\label{orbit}. \emph{Given  $\ax$,~ $\la$,~ $\px$ as above,
then
\\
1) ~ the following formula for the character $\chi_\la$ of
representation $M(\la)$ holds
\begin{equation}\label{ww}\chi_\la(1+x) = \frac{1}{\sqrt{|\Omega|}}
\sum_{\mu\in\Omega(\la)} e^{\mu(x)},\quad\quad
x\in\ax;\end{equation} 2)~
$\dim M(\la) = q^{\codim\,\px} = \sqrt{|\Omega|}$;\\
3) the representation   $M(\la)$ doesn't depend on a choice  of associative polarization;\\
4) the representation  $M(\la)$ is irreducible;\\
5) let two linear forms  $\la$ and $\la'$ have  associative
polarizations;  the representations  $M(\la)$ and $M(\la')$ are
equivalent if and only if  $\la$ and $\la'$ lie in a common
coadjoint orbit.}\\
\Proof.  For any   $\la\in\ax^*$, we denote by  $\ax^\la =
\{a\in\ax:~~ \la([a,\ax]) = 0\}$ the stabilizer of $\la$ in the Lie
algebra $\ax$. Obviously, the equality  $\la((E+a)x)=\la(x(E+a))$
 is equivalent to   $\la(ax)=\la(xa)$. Hence, $E+\ax^\la$ coincides
 with the stabilizer  $G^\la$  of linear  form  $\la$ in the group  $G$.
 It implies
 \begin{equation}\label{OO}|\Omega| = \frac{|G|}{|G^\la|} =
\frac{|\gx|}{|\gx^\la|} = q^{\dim \gx - \dim \gx^\la} =
q^{\dim\Omega}\end{equation} The proof may be finished similarly as
in papers \cite{A2, Kzh, Sa, P1}. $\Box$
\\
\Remark. The main result of the paper  \cite{IK} implies that  the
formula (\ref{ww}) in general does not true for  $\ut(n,\Fbq)$.
Therefore, it is not true that every linear form  on $\ut(n,\Fbq)$
has an associative polarization. Applying the classification of
coadjoint orbits for  the unitriangular Lie groups of lower sizes
\cite{P2}, one can prove the existence of  associative polarization
for $n\leq 7$.

\section{\bf {Orbits and representations for  $\laSa$}}
Let $\la$ be an nondegenerate character of  $\gx_+$  as an
associative algebra. In this section, we construct  families of the
linear forms $\laSa$ where $S$ is a subset of $\Pi$ and $a\in
\La_S$. For any $\laSa$, we construct  an associative polarization
$\px_S$  and corresponding irreducible representation $\VaSa$. In
what follows, we obtain a description of the coadjoint orbit $\OSa$;
we shall show that the representations $\{\VaSa\} $ are pairwise
non\-equiv\-alent.

  We shall treat any root $\gamma = (i,j)$ as a box in empty   $n\times n$-matrix.
  The number  $i$ is called the row of  the root
$\gamma$, and $j$, respectively, the column of  the root $\gamma$.
We shall say that a root  $\gamma'=(i',j')$ lies on the left side
(respectively, stronger on the left side) of the root $\gamma$, if
$j'\leq j$ (respectively, $j'<j)$. We define similarly relations of
lying on the right side, over and lower.

Let  $S$ be an arbitrary subset of  $\Pi$. \\
\Notations.  \\
1) ~ Denote by   $L_S^0$ a subset that consists of all roots $\gamma
=(i,j)$ obeying the following conditions:\\
i)~ $1\leq i\leq k$ and $i+j\geq n+1$,\\
ii)~ there is no roots of  $S$ in the $i$th row and $j$th column,\\
iii)~ all roots of  $\Pi$ that are lying stronger over and stronger
on the left side of $\gamma$  belong to  $S$.\\
2) The subset  $L_S^{00}$  is empty for  the odd  $n$. For
$n=2(k+1)$, the subset  $L_S^{00}$ consists of the single root
$\gamma=(k+1,j)$ obeying conditions ii) and  iii). \\
3)  ~ $L_S^+$. A root $\gamma$ of $ R_0\sqcup R_-$ belongs to
$L_S^+$ , if it  belongs to  the same column
and lies stronger over   some root of  $L_S^0$.\\
4)   ~ $L_S^-$. A root $\gamma'=(n-i,j)$ belongs to  $L_S^-$ if and
only if $\gamma =(i,j) $ belongs to  $L_S^+$. \\
5) ~$L_S = L_S^+\sqcup L_S^0\sqcup L_S^{00}\sqcup L_S^-$.\\
6)~ $R_S = R_+\sqcup L_S$;\\
7)~ $|S|=s$,~ $|R_+|=r_+$,~ $|R_0|=r_0$.

Notice that    $|L_S^+| = |L_S^-| = |S| = s$, ~$|L_S^0| = k-s$, ~
$|L_S^{00}| = \eps$,~ $k+\eps = r_0$

Consider the following ordering in the set all roots $R$: we say
that $\beta\geq \al$, if $\beta$ lies lower than  $\al$, or in the
same row and on the left side from  $\al$.
 Order the subset  $R_+\setminus S$ with respect to this  ordering
$R_+\setminus S=\{\al_1<\al_2<\ldots<\al_{r_+-s}\}$.\\
\Lemma\Num\label{ll}. \\\emph{ 1) For any root $\al_i$, there exists
a unique root  $\beta _i\in R\setminus R_S$ such that
$\al_i+\beta_i$ belongs to $
S\sqcup L_S^0$.\\
2) ~ For any $1\leq i < j\leq r_+-s$,  the sum  $\al_i+\beta_j$
 is either undefined, or is defined and  belongs to  $\left(R_+\setminus \{S\sqcup \Pi_0\}\right)\sqcup
L_S^+$.}\\
\Proof~ Follows from the definitions. $\Box$

Denote by $\lx_S$ (resp. $\lx_S^\pm$, $\lx_S^0$, $\lx_S^{00}$) a
subspace spanned  by the system  $E_\gamma$,~ $\gamma\in L_S$ (resp.
$\gamma \in L_S^\pm$, $\al\in L_S^0$). Obviously,
$$\lx_S=\lx_S^+\oplus \lx_S^0\oplus \lx_S^{00}\oplus \lx_S^-.$$ The
subspace  $\lx_S$ is an associative subalgebra with zero
multiplication. The subspace
$$\px_S = \spp\{E_\gamma:~ \gamma\in
 R_S\} = \gx_+\oplus\lx_S$$ is also an associative subalgebra in  $\gx$,
and  $P_S = E + \px_S$ is a subgroup in  $G$.

Consider the subset  $\La_S$ that consists of all functions
$$a: L_S^0\sqcup L_S^{00}\sqcup L_S^-\to \Fbq$$
such that   $a(\gamma)\ne 0$ for any $\gamma\in L_S^0$. One may
identify
 $$\La_S \cong \Fbq^{s+\eps}\times\left(\Fbq^*\right)^{k-s}.$$
The number of elements of  $\La_S$ equals to
$q^{s+\eps}(q-1)^{k-s}$.

For any  $a\in\La_S$, we define a linear form  $\la_{S,a}$ on $\gx$
as follows:\\
1) ~ $\laSa(E_\gamma) = \la(E_\gamma)$ for all  $\gamma\in
\Pi_0\setminus\Pi$ and $\gamma\in
S$,\\
2)~ $\laSa(E_{\gamma}) = a(\gamma)$ for all  $\gamma\in
L_S^0\sqcup L_S^{00}\sqcup L_S^-$,\\
3) ~ $\laSa(E_{\gamma}) = 0$ for all other  $E_\gamma\in\gx$.

Notice that the definition of  $\La_S $ implies \\
i) ~ $\laSa(E_\gamma) \ne 0$ for all  $\gamma\in L_S^0$;\\
ii)~  $\laSa(E_\gamma)$ may have an arbitrary values, when
$\gamma\in L_S^{00}\sqcup L_S^{-}$.\\
iii) ~ $\laSa(E_\gamma)  = 0$ for all  $\gamma$ of $ R_+\setminus
\{S\sqcup \Pi_0\}$ and all $\gamma\in L_S^+$.

Easy to see that the linear form  $\laSa$ is a character of the
associative algebra $\px_S$ in the sense of definition  \ref{chi}.
Following formula (\ref{xi}), we define a complex character
$\xi_{\la,S,a}$ of the subgroup  $P_S=E+\px_S$. We denote
$$V_{S,a}(\la) = \ind(\xi_{\la,S,a},P_S,G).$$
\Theorem\Num\label{th1}. \emph{Let  $\la$ be a nondegenerate
character of the algebra $\gx_+$ (see definition \ref{nondeg}). Then\\
1) the subalgebra  $\px_S$ is an associative polarization for the
linear form  $\laSa$,\\
2)  every representation  $V_{S,a}(\la) $ is irreducible.} \\
\Proof. By proposition  \ref{orbit}, the statement  1) implies 2)
1). Let us prove 1).

The subspace $\px_S$ is an associative subalgebra,  and
$\laSa(\px_S^2)=0$. It is sufficient to prove that  $\px_S$ is a
maximal isotropic subspace for the skew symmetric bilinear form
$\laSa([x,y])$. Suppose the contrary. Assume that  there exists  $ x
\in \gx\setminus \px_S$ such that  $\laSa([\px_S, x])=0$. Then
$$x=\sum_{j=1}^{r_+-s} b_jE_{\beta_j}.$$
Let  $i$ be the smallest number obeying  $b_i\ne 0$. Lemma \ref{ll}
and definition of $\laSa$ imply that
$$\laSa([E_{\al_i},E_{\beta_i}]) = c_i\ne 0,\quad
\laSa([E_{\al_i},E_{\beta_j}]) = 0$$ for all $j>i$.

Hence $\laSa([\px_S, x])=c_i b_i\ne 0$. This leads to contradiction.
Therefore,  $\px_S$ is a maximal isotropic subspace. $\Box$

Denote by  $\OSa$ the orbit of   $\laSa\in\gx^*$ with respect to the
coadjoint  representation of  the group  $G$. \\
Introduce the following notations: ~ $\oFbq$ is an algebraic closure
of the field $\Fbq$, ~ $\ogx = \gx\otimes \oFbq$, ~ $\oOSa$ is the
coadjoint orbit of  $\laSa$ with respect to the group $\overline{G}
= \UT(n,\oFbq)$. The orbit  $\oOSa$ is closed, since every orbit of
a regular action of a nilpotent group on an arbitrary affine
algebraic variety is closed  (\cite{Dix}, 11.2.4). We shall find a
system of generators of the defining ideal $\IOSa$ of the orbit
$\oOSa$.

Notice that  $\dim \oOSa = 2 \codim\, \px_S =2(r_+-s)$.
Respectively,
$$\codim\,\oOSa = \dim \gx - \dim \oOSa = 2r_++r_0-2(r_+-s)=r_0+2s.$$
Notice that the number of roots in  $S\sqcup L_S$ also equals to
$r_0+2s$, since  $|S|=|L_S^+|= |L_S^-|=s$, and  $|L_S^0| = r_0-s$.
In what follows, we shall correspond some element of the symmetric
algebra $\Sc(\ogx)=\oFbq[\ogx^*]$ to  any  root  $\gamma\in S\sqcup
L_S$, and we shall prove  that the constructed system of elements
will generate the defining ideal $\IOSa$.

Denote by   $\Xb$ the upper triangular matrix with zeros on the
diagonal and the following entries over the diagonal:  any place
$(i,j)$, where $1\leq i<j\leq n$, is filled by the matrix unit
$E_{ij}$. Every minor of the matrix $\Xb$ is an element of the
symmetric algebra $\Sc(\ogx)$, that is a polynomial on $\ogx^*$.
Given a root $\gamma\in R$, we consider the system $\Sb_\gamma$ that
consists of the root $\gamma$ and also  of all roots from $ S\sqcup
L_S^0$ lying strongly over and strongly on the right side from
$\gamma$.  Denote by $\Mc_\gamma$ the minor of the matrix $\Xb$ that
has the systems of rows and columns just as $\Sb_\gamma$ has.

 Let $\gamma=(n-i,j)\in L_S^-$. Consider the characteristic matrix $\Xb-\tau
 E$; cutting first $i$ columns  and last $i$ rows, we obtain the minor
   $\left|\Xb-\tau E\right|_{i}$.  Then
 $$ \left|\Xb-\tau E\right|_{i} = \Pc_{\gamma,0} \tau^{n-2i} + \Pc_{\gamma,1}
 \tau^{n-2i-1}+\ldots+\Pc_{\gamma,n-2i}.$$
  Denote   $$\Fc_\gamma = \left\{\begin{array}{l} \Mc_\gamma,
 ~~\mbox{if}
 ~~ \gamma\in S\sqcup L_S^+\sqcup L_S^0\sqcup L_S^{00}, \\
 \Pc_{\gamma,1}, ~~\mbox{if}
 ~~ \gamma\in  L_S^-.
 \end{array}\right.$$
Denote by $\Fc_\gamma^0$  a value of  $\Fc_\gamma$ at the point  $\laSa$. \\
\Remark\Num\label{rr}. Notice that  $\Fc_\gamma^0 =
c\laSa(E_\gamma)$. Here  $c$ equals to  a  product of values of
$\laSa$ on some    $E_\mu$, where $\mu<\la$ and $\mu$ belongs to
$S\sqcup
L_S^0$. Therefore, $c\ne 0$. \\
 \Theorem\Num\label{th2}. \emph{The defining ideal  $\IOSa$ of the orbit
$\oOSa$ is generated by  the algebraically independent system of
polynomials
\begin{equation}\label{ge} \{\Fc_\gamma - \Fc_\gamma^0:~~ \gamma\in S\sqcup
L_S\}\end{equation}}\\
 \Proof.
As above we order   $R  = \{\gamma_1<\ldots<\gamma_N\}$, where
$N=\frac{n(n-1)}{2}$,  with respect to the ordering introduced above
(before lemma  \ref{ll}).  The ordering on the set of $R$ provides
the ordering on the set of all matrix units $$\{E_{i,j}:~ 1\leq
i<j\leq n\}.$$ The associative algebra $\gx$ has a chain of ideals
$$\ogx_1= <E_{1,n}>\subset\ogx_2\subset\ldots\subset\ogx_N=\gx,$$ where
 $\gx_i$ is a span of all matrix units with numbers  $\leq i$.
  Denote by  $\Ic_i$ the ideal in $\Sc(\ogx)$ generated by  whose
  elements of   (\ref{ge}) that has number $\leq i$. The last ideal $\Ic_N$ coincides with
  the ideal $\Ic$ generated by  the system of generators  (\ref{ge}).

  It is not difficult to show that
   \begin{equation}\label{ff}
  \Fc_{\gamma_i} = cE_{\gamma_i} + \Phi_{i-1}\bmod \Ic_{i-1},\end{equation}
where $c$ is the  constant as in remark above, $\Phi_{i-1}$ is some
polynomial of
  $S(\ogx_{i-1})$. Any  $E_{\gamma_i}$ belongs to  $\ogx_i$, does not belong to
   $\ogx_{i-1}$; using  (\ref{ff}), we conclude that the system of generators
   (\ref{ge}) is algebraically independent, and the ideal $\Ic$ is prime.
   The number of generators of the set (\ref{ge}) equals to  $r_0+2s$; this is equal to codimension
    of the orbit  $\OSa$. Therefore,
  $\dim\,\mathrm{Ann}\,\Ic = \dim\OSa$.

  Obviously,  the generators (\ref{ge}) annihilate at the point
   $\laSa$. To finish the proof, it is sufficient to show that
    the ideal  $\Ic$ is invariant with respect to the adjoint representation of the group
    $G$. Using  direct calculations, one can show that  for any  $1\leq m\leq n-1$
the element   $(\mathrm{ad}\,E_{m,m+1})\Fc_{\gamma_i} $ belongs to
the ideal  $\Ic_{i-1}$.  $\Box$\\
 \Theorem\Num\label{th3}. \emph{Linear forms  $\laSa$ и $\laSas$ lie in a common
 $\Ad_G^*$-orbit if and only if they coincide.}\\
\Proof.  It is obvious that, if $\laSa=\laSas$, then they are lying
in the common orbit. Let us prove the contrary.

Suppose that  $\laSa$ and $\laSas$  are lying in a common
$\Ad_G^*$-orbit. Then  they are lying in a common $\Ad^*$-orbit with
respect to the group  $\overline{G}$. The defining ideal of the
common orbit is generated by  the system of polynomials (\ref{ge}).
The values of every polynomial of  $\{\Fc_\gamma\}$ at the points
$\laSa$ and $\laSas$ coincide.

Suppose that  $S\ne S'$. Order the set  $\Pi$ in accordance with the
number of row. Choose a number  $i$ such that  \\
i)~ the subset of all roots  of  $S$ with number of row  $< i$
coincides with the same subset of  $S'$;\\
ii)~ the root  $\gamma=(i,n-i)$ belongs to  $S$ and does not belong
to  $S'$.

Since  $\gamma\in S$, there is some  $\gamma_* =(i,j_*)\in L_S^+$
that   lies strongly  on the right side from $\gamma$. In accordance
with remark  \ref{rr}, a value  $\Fc_{\gamma_*}^0$ of the polynomial
$\Fc_{\gamma_*}$  at the point  $\laSa$ equals to  $c
\laSa(E_{\gamma_*})$, where $c\ne 0$. Since $\gamma_*\in L_S^+$, we
have $\laSa(E_{\gamma_*})=0$ and, therefore,  $\Fc_{\gamma_*}^0 =0$.

 On the other hand, since   $\gamma\notin S'$, we have  $\gamma_* \in L_{S'}^0$.
Arguing similarly, we obtain  $\Fc_{\gamma_*}^0 \ne 0$. This leads
to contradiction. Hence,  $S=S'$.

The equality   $S=S'$ implies  $L^\pm_S=L_{S'}^\pm$ and
$L^0_S=L_{S'}^0$. Calculating values $\Fc_\gamma^0$ for $\gamma\in
L_S^0$, we conclude that  $a=a'$. $\Box$

\section{\bf {Decomposition into irreducible components}}

\Theorem\Num\label{th4}.\\
\emph{1) Every representation  $V_{S,a}(\la)$ can be realized as a
subrepresentation of  the representation   $V(\la)$. \\
2) Any representation of Gelfand-Graev type  is decomposed into a
sum of irreducible components
\begin{equation} \label{vla} V(\la) = \bigoplus \VaSa,\end{equation} where
$S\subset \Pi$ and $a\in \La_S$. The multiplicity of every
irreducible component of $V(\la)$ is equal to one.
} \\
\Proof. Let   $v$ be an eigenvector in  $\VaSa$ for $G_+ = \{ 1+x:~
x\in\gx_+\}$; its eigenvalue has a form  $e^{\nu(x)}$, where
$\nu(x)$  is a character of the associative algebra $\gx_+$. Then
$\nu(x)$ is called an $\gx_+$-weight, and  corresponding eigenvector
$v$ is called an $\gx_+$-weight vector.

To prove statement  1), it is sufficient to show that there exists
nontrivial homomorphism $V(\la)\to\VaSa$. This is equivalent to
existence of nonzero $\gx_+$-weight vector  of weight  $\la$ in
$\VaSa$. The proof of this statement is similar to proof of
proposition  \ref{VV}. Let $f_0$ be the generating vector of
representation  $\VaSa$. For any root $\al\in \Pi\setminus S$, there
exists a unique root $\gamma(\al)=(i,j_*)\in L_S^0$ that lies
strongly on the right side from  $\al$. Then  $\gamma(\al) =
\al+\beta(\al)$, where $\beta = (n-i,j_*)\in R_-$. The subspace,
spanned by the set  $E_{\beta(\al)},~~\al\in \Pi\setminus S$, is a
subalgebra with zero multiplication. Hence, for all  $\al,\al'\in
\Pi_0\setminus \Pi$, the elements of one-parameter subgroups
$x_{\beta(\al)}(t)$ and $x_{\beta(\al')}(t')$ commute. Recall that
$\laSa(E_\al)=0$ and $\laSa(E_{\gamma(\al)})=a(\gamma(\al))\ne 0$.
 Consider the element
\begin{equation} g_0 = \prod x_{\beta(\al)}(t_\al)\in
G,\end{equation} where  $\al$ is running through $\Pi\setminus S$,
and $t_\al\in \Fbq$ is a solution of the equation
$$ \laSa(E_{\gamma(\al)})t_\al = \la(E_\al).$$ The vector  $g_0f_0$ is a
$\gx_+$-weighted vector with the weight $\la$. This proves  1).

Taking into account theorems \ref{th1}, \ref{th2}, to prove the
statement  2), it is sufficient  to show that dimensions of
representations  of  left and right hand sides of formula
(\ref{vla}) coincide. Denote
\begin{equation}\label{nnn}
V_S(\la) = \bigoplus_{a\in\La_S} \VaSa. \end{equation} We  have как
$|\La_S|= q^{s+\eps}(q-1)^{k-s}$, ~~ $\dim \VaSa = q^{r_+-s}$, and
$$\dim V_S(\la) = q^{r_++\eps}(q-1)^{k-s}.$$ Then
$$\dim\left(\bigoplus_{S\subset \Pi} V_S(\la)\right) =
q^{r_++\eps}\sum_{s=0}^k C_k^s (q-1)^{k-s} =
q^{r_++k+\eps}=q^{r_++r_0} = \dim V(\la).$$ $\Box$\\
\Cor\Num\label{mmm}.\emph{ The number of irreducible components in
the representations of Gelfand-Graev type  $V(\la)$ for the group
$\UT(n,\Fbq)$  does not depend on choice of nondegenerate character
 $\la$, and it equals to $q^\eps(2q-1)^k$.\\}
\Proof. The number of irreducible components of  $V(\la)$ is equal
to
$$\sum_{s=0}^kC_s^k |\La_S| = \sum_{s=0}^k
C_k^sq^{s+\eps}(q-1)^{k-s} = q^\eps \sum_{s=0}^k
C_k^sq^{s}(q-1)^{k-s} = q^\eps(2q-1)^k.$$

\section{\bf {The Hecke algebra}}

Let  $G$ be an arbitrary finite group, $\Ac_G$ be its group algebra
over  $\Cb$, ~ $H$ be a subgroup of $G$, ~$\xi$ be a character
(one-dimensional representation) of the group $H$. Denote by $P_\xi$
the element $\sum_{h\in H}\xi(h^{-1})h$ in the group algebra
$\Ac_G$.

Let us realize the induced representation  $V=\ind(\xi,H,G)$ in the
space   $\Ac_GP_\xi$ by left multiplication. The Hecke  algebra
$\Hc(V)$ (i.e. the the algebra of all $G$-endomorphisms of $V$) is
isomorphic to subalgebra $P_\xi\Ac_GP_\xi$ with inverted
multiplication. The algebra $P_\xi\Ac_GP_\xi$ is spanned by the
system  $\{P_\xi xP_\xi:
~x\in G\}$. The following is well known  (see \cite[Lemma 84]{S}).\\
1) ~ The element  $P_\xi xP_\xi$ is uniquely determined up to
constant nonzero multiple determined be double class $HxH$. So the
algebra $\Hc(V)$ is spanned by the elements  $P_\xi xP_\xi$, where
$x$ is running
through a system of representatives of double classes $HxH$. \\
2) ~ The elements  $P_\xi xP_\xi$, where $x$ is running through a
system of whose representations of double  $(H,H)$ classes that
satisfy  $P_\xi xP_\xi\ne 0$, form a basis of  $\Hc(V)$.\\
3)  Define a character  $x\xi$ on $ xHx^{-1}$ by the formula
$x\xi(y)=\xi(x^{-1}yx)$. The element  $P_\xi xP_\xi$ is nonzero if
and only if $\xi = x\xi$ on $xHx^{-1}\cap H$. Summarizing  1)-3), we
obtain
the following proposition.\\
\Prop\Num\label{xl}. \emph{Let   $V=\ind(\xi,H,G)$. Then the system
of elements $P_\xi xP_\xi$, where $x$ is running through a system of
whose representatives of the double  $(H,H)$ classes that satisfy
$\xi = x\xi$ on $xHx^{-1}\cap H$, form a basis of  $\Hc(V)$.}

Turn to representations of Gelfand-Graev type. We denote by
$\Hc(\la)$ the Hecke algebra of the representation $V(\la)$. Since
$V(\la)$ is multiplicity free, the Hecke is commutative. The
dimension of Hecke
 algebra equals to the number of irreducible components of
$V(\la)$, it is equal to  $q^\eps (2q-1)^k$. Put $H=G_+$ and $\xi
=\xi_\la$ (see formula (\ref{chi})). Our goal is to construct a
system of elements $\{x\}$ in $G$ such that  $P_\xi x P_\xi$ is a
basis in  $\Hc(\la)$.

 Given a subset
$S\subset \Pi$ , we define  the subset  $\La_S'$ that consists of
all vectors $(b_1,\ldots,b_{k+\eps})$, where all $b_i\in \Fbq$, and
$b_i\ne 0$, if $(i,n-i)\in S$.

We construct the matrix  $X_{S,b}= (x_{ij})$ as follows. \\
1)~ The matrix $X_{S,b}$  lies in $ E + \gx_0 + \gx_-$, i.e.
$x_{ii}=1$ and $x_{ij}=0$  for all  $i>j$ and for whose pairs   $i<j$ that obey   $i+j<n+1$. \\
2) First, we fill the last column of the matrix  $X_{S,b}$. For any
$1\leq i\leq k+\eps$ we put   $x_{in}= b_i$. There is a unique root
$(i,n-i)\in\Pi$ in each row $1\leq i\leq k$.
 If  $(i,n-i)\notin
S$, then we put  $x_{n-i,n}=0$. If $(i,n-i)\in S$, then we put
$x_{n-i,n}=b_i$. So $x_{i,n}=x_{n-i,n}$, if $1\leq i\leq
k$ and $(i,n-i)\in S$. \\
3)~ Now, we fill the other columns. If $(i,n-i)\in S$, then we put
$x_{i+1,n-i} =\ldots = x_{n-i-1,n-i}=0$.  If $(i,n-s)\notin S$, then
we put
$$\left(\begin{array}{c}
  x_{i+1,n-i}  \\
         \vdots \\
          x_{n-i-1,n-i}\\
       \end{array}\right) = x_{i,n}
\left(\begin{array}{c}
  x_{i+1,n-i}  \\
         \vdots \\
          x_{n-i-1,n-i}\\
       \end{array} \right).$$
\Theorem\Num\label{th5}. \emph{The system of elements
\begin{equation}\label{PP}\{P_\xi X_{S,b}P_\xi:~ S\subset
\Pi, b\in \La'_S\}\end{equation}  is a basis of  $\Hc(V)$.}\\
\Proof. By direct calculations, we show that any  $x\in \{X_{S,b}\}$
obeys  $\xi = x\xi$ on $xHx^{-1}\cap H$, where $H=G_+$ and
$\xi=\xi_\la$. After, we  prove that  elements of $\{X_{S,b}\}$ lie
in different double classes. The number of elements of (\ref{PP}) is
equal to
$$\sum_{s=0}^k C_k^s(q-1)^sq^{k+\eps} = q^\eps (2q-1)^k, $$ i.e. it is equal to  $\dim\Hc(V)$.
By proposition \ref{xl}, this concludes the proof. $\Box$

\section{\bf {Calculations for small  $n$}}

In all examples below, using the Killing form,  we identify $\gx^*$
with the space of lower triangular matrices with zeros on the
diagonal.\\
\Ex\Num. Case  $n=3$. In this case,
 $$\gx = \left\{\left(%
\begin{array}{ccc}
  0 & x_{12} & x_{13} \\
  0 & 0 & x_{23} \\
  0 & 0 & 0 \\
\end{array}%
\right)\right\}, \quad  R_+=\{(1,2)\}, \quad  R_0=\{(1,3)\}, \quad
R_-=\{(2,3)\},$$
$$
\gx_+ = \left\{\left(%
\begin{array}{ccc}
  0 & x_{12} &  0 \\
  0 & 0 & 0 \\
  0 & 0 & 0 \\
\end{array}%
\right) \right\}, \quad
\gx_0 = \left\{\left(%
\begin{array}{ccc}
  0 & 0 & x_{13} \\
  0 & 0 & 0 \\
  0 & 0 & 0 \\
\end{array}%
\right)\right\},\quad
\gx_- = \left\{\left(%
\begin{array}{ccc}
  0 &0 & 0 \\
  0 & 0 & x_{23} \\
  0 & 0 & 0 \\
\end{array}%
\right)\right\},
$$
The set   $\Pi$ coincides with   $R_+ = \{(1,2)\}$. Nondegenerate
character  $\la$ on $\gx_+$ is defined by one number
$\la(E_{12})=c\ne 0$. There exist  only two subsets in $\Pi$: the
empty set and $\Pi$.
Recall that, in the odd case,  $L_S^{00}=\varnothing$.\\
i)~ $S = \varnothing $. Then  $L_S^+ = L_S^-=\varnothing$,
$L_S^0=\{(1,3)\}$, ~ $\La_S \cong\Fbq^* = \{a\in\Fbq:~ a\ne 0\}$,

$$\px_S = \left\{\left(%
\begin{array}{ccc}
  0 & x_{12} & x_{13} \\
  0 & 0 & 0 \\
  0 & 0 & 0 \\
\end{array}%
\right) \right\}, \quad \la_{S,a} =
\left(%
\begin{array}{ccc}
  0 & 0 &  0 \\
  0 & 0 & 0 \\
  a & 0 & 0 \\
\end{array}%
\right)
$$
The coadjoint orbit  $\OSa$ is defined in  $$\gx^* = \left\{\left(%
\begin{array}{ccc}
  0 & 0 & 0 \\
  y_{21} & 0 & 0 \\
  y_{31} & y_{32} & 0 \\
\end{array}%
\right) \right\}~~~\mbox{by ~equation}~~ y_{31}=a.$$ \\
ii) $S = \{(1,2)\}$. Then ~~$L_S^+ = \{(1,3)\}$,~~ $ L_S^-=
\{(2,3)\}$, ~~ $L_S^0=\varnothing$, \\ $\La_S \cong\Fbq =
\{a\in\Fbq\}$,

$$\px_S = \left\{\left(%
\begin{array}{ccc}
  0 & x_{12} & x_{13} \\
  0 & 0 & x_{23} \\
  0 & 0 & 0 \\
\end{array}%
\right) \right\}, \quad \la_{S,a} =
\left(%
\begin{array}{ccc}
  0 & 0 &  0 \\
  c & 0 & 0 \\
  0 & a & 0 \\
\end{array}%
\right)
$$

The coadjoint orbit  $\OSa$ is defined in $\gx^*$ by the equations \\
$y_{31}=0$,~ $y_{21}=c$,~ $y_{32}=a$. Any representation of
Gelfand-Graev type  $V(\la)$ decomposes into a sum of irreducible
representations that correspond to  mentioned above orbits. The
number of irreducible components equals to  $2q-1$. The matrices
$X_{S,b}$ from theorem  \ref{th5} has the  form
$$
\left(%
\begin{array}{ccc}
  1 & 0 & b \\
  0 & 1 & b \\
  0 & 0 & 1 \\
\end{array}%
\right), ~~\mbox{where}~~ b\in\Fbq^*,\quad\quad
\left(%
\begin{array}{ccc}
  1 & 0 & a \\
  0 & 1 & 0\\
  0 & 0 & 1 \\
\end{array}%
\right), ~~\mbox{where}~~ a\in\Fbq.$$ \Ex\Num. Case $n=4$. In this
case, $ R_+=\{(1,2), (1,3)\}$,~$ R_0=\{(2,3),(1,4)\}$,~ \\$
R_-=\{(2,4), (3,4)$,
  $$\gx = \left\{\left(%
\begin{array}{cccc}
  0 & x_{12} & x_{13} &x_{14}\\
  0 & 0 & x_{23} &x_{24} \\
  0 & 0 & 0 &x_{34}\\
0 & 0 & 0 & 0\\
\end{array}%
\right)\right\},  \quad
\gx_+ =  \left\{\left(%
\begin{array}{cccc}
  0 & x_{12} & x_{13} & 0\\
  0 & 0 & 0 &0 \\
  0 & 0 & 0 &0\\
0 & 0 & 0 & 0\\
\end{array}%
\right)\right\},$$ $$
\gx_0 = \left\{\left(%
\begin{array}{cccc}
  0 & 0 & 0 &x_{14}\\
  0 & 0 & x_{23} & 0 \\
  0 & 0 & 0 & 0\\
0 & 0 & 0 & 0\\
\end{array}%
\right)\right\},\quad
\gx_-= \left\{\left(%
\begin{array}{cccc}
  0 & 0 & 0 & 0\\
  0 & 0 & 0 &x_{24} \\
  0 & 0 & 0 &x_{34}\\
0 & 0 & 0 & 0\\
\end{array}%
\right)\right\}
$$
The set   $\Pi$ coincides with  $ \{(1,3)\}$. A nondegenerate
character  $\la$ on $\gx_+$ is defined by one number
$\la(E_{13})=c\ne 0$. There exist only two subsets in  $\Pi$: the
empty
subset and  $\Pi$.\\
i)~ $S = \varnothing $. Then  $L_S^+ = L_S^-=\varnothing$, $L_S^0=
\{(1,4)\}$, ~  $L_S^{00}=\{(2,3)\}$,~ $\La_S \cong\Fbq^* \times \Fbq
= \{(a_1,a_2)\in\Fbq^2:~ a_1\ne 0\}$,
$$\px_S = \left\{\left(%
\begin{array}{cccc}
  0 & x_{12} & x_{13} & x_{14}\\
  0 & 0 & x_{23}&0 \\
  0 & 0 & 0 &0\\
0 & 0 & 0 & 0\\
\end{array}%
\right)\right\}, \quad \la_{S,a} =
\left(%
\begin{array}{cccc}
  0& 0 & 0 & 0 \\
  0 & 0 & 0 & 0 \\
  0 & a_2 & 0 & 0 \\
  a_1 & 0 & 0 & 0 \\
\end{array}%
\right). $$
The coadjoint orbit  $\OSa$ if defined in  $$\gx^* = \left\{\left(%
\begin{array}{cccc}
  0 & 0 & 0 &0\\
  y_{21} & 0 & 0&0 \\
  y_{31} & y_{32} & 0&0\\
y_{41} & y_{42} & y_{43}&0   \\
\end{array}%
\right) \right\}~~~\mbox{by~ the ~equations }~~ y_{41}=a_1,~~\left|
\begin{array}{cc}
  y_{31} & y_{32} \\
  y_{41} & y_{42} \\
\end{array}
\right| = -a_1a_2.$$
 ii) $S = \{(1,3)\}$. Then ~~$L_S^+ = \{(1,4)\}$,
~~ $L_S^-=\{(3,4)\}$,
 ~ $L_S^0=\varnothing$,~\\ $L_S^{00}=\{(2,4)\}$,~
$\La_S \cong  \Fbq^2 = \{(a_1,a_2)\in\Fbq^2\}$,
$$\px_S = \left\{\left(%
\begin{array}{cccc}
  0 & x_{12} & x_{13} & x_{14}\\
  0 & 0 & 0&x_{24} \\
  0 & 0 & 0 &x_{34}\\
0 & 0 & 0 & 0\\
\end{array}%
\right)\right\}, \quad \la_{S,a} =
\left(%
\begin{array}{cccc}
  0& 0 & 0 & 0 \\
  0 & 0 & 0 & 0 \\
  c & 0 & 0 & 0 \\
  0 & a_1 & a_2 & 0 \\
\end{array}%
\right). $$ The coadjoint orbit  $\OSa$ is defined in  $\gx^*$ by
the equations $y_{41}=0$, ~$y_{42}=a_1$, ~$y_{31}=c$, ~ $
y_{42}y_{21} + y_{43}y_{31} = a_2c $.

The coadjoint orbit  $\OSa$ is defined in $\gx^*$ by the equations \\
$y_{31}=0$,~ $y_{21}=c$,~ $y_{32}=a$. Any representation of
Gelfand-Graev type $V(\la)$ decomposes into a sum of irreducible
representations that correspond to  mentioned above orbits. The
number of irreducible components equals to  $q(2q-1)$. The matrices
$X_{S,b}$ from theorem  \ref{th5} has the  form
$$
\left(%
\begin{array}{cccc}
  1 & 0 &0& b \\
  0 & 1 &0& a \\
  0 & 0 &1& b\\
0 & 0 & 0&1 \\
\end{array}%
\right), ~~\mbox{where}~~ b\in\Fbq^*,~ a\in\Fbq, \quad  \left(%
\begin{array}{cccc}
  1 & 0 &0& b \\
  0 & 1 &ba&a  \\
  0 & 0 &1& 0\\
0 & 0 & 0&1 \\
\end{array}%
\right), ~~\mbox{where}~~  a,b\in\Fbq.$$ \Ex\Num. Case  $n=5$. In
that case  $ R_+=\{(1,2), (1,3), (1,4), (2,3)\}$,~ \\$
R_0=\{(2,4),(1,5)\}$,~$ R_-=\{(2,4), (3,4), (3,5), (4,5)\}$,
{\small $$\gx = \left\{\left(%
\begin{array}{ccccc}
  0 & x_{12} & x_{13} &x_{14}&x_{15}\\
  0 & 0 & x_{23} &x_{24}&x_{25} \\
  0 & 0 & 0 &x_{34}&x_{35}\\
0 & 0 & 0 & 0&x_{45}\\
0 & 0 & 0 & 0&0\\
\end{array}%
\right)\right\},  \quad
\gx_+ =  \left\{\left(%
\begin{array}{ccccc}
  0 & x_{12} & x_{13}&x_{14}& 0\\
  0 & 0 &  x_{13}&0&0 \\
  0 & 0 & 0 &0&0\\
0 & 0 & 0 & 0&0\\
0&0&0&0&0\\
\end{array}%
\right)\right\},$$ $$
\gx_0 = \left\{\left(%
\begin{array}{ccccc}
  0 & 0 & 0 &0&x_{15}\\
  0 & 0 &0& x_{24} & 0 \\
  0 & 0 & 0 & 0&0\\
0 & 0 & 0 & 0&0\\
0&0&0&0&0\\
\end{array}%
\right)\right\},\quad
\gx_-= \left\{\left(%
\begin{array}{ccccc}
  0 & 0 & 0 & 0&0\\
  0 & 0 & 0 &0&x_{25} \\
  0 & 0 & 0 &x_{35}&x_{35}\\
0 & 0 & 0 &0& x_{45}\\
0&0&0&0&0\\
\end{array}%
\right)\right\}
$$}
The set $\Pi$ coincides with   $\{(2,3), ~(1,4)\}$. A nondegenerate
character  $\la$ on  $\gx_+$ is defined by pair of numbers
$\la(E_{14})=c_1\ne 0$ и $\la(E_{23})=c_2 $ (the second number may
by   arbitrary). There are  four subsets in  $\Pi$.\\
i)~ $S = \varnothing $. Then  $L_S^+ = L_S^- = \varnothing$,~~
$L_S^0= \{(1,5), (2,4)\}$,~~ $\La_S \cong\Fbq^{*2} =\\
\{(a_1,a_2)\in\Fbq^2:~ a_1, a_2\ne 0\}$,
{\small$$\px_S =  \left\{\left(%
\begin{array}{ccccc}
  0 & x_{12} & x_{13} & x_{14}& x_{15}\\
  0 & 0 & x_{23}& x_{24}&0 \\
  0 & 0 & 0 &0&0\\
0 & 0 & 0 & 0&0\\
0&0&0&0&0\\ \end{array}%
\right)\right\}, \quad \la_{S,a} =
\left(%
\begin{array}{ccccc}
0&  0& 0 & 0 & 0 \\
 0& 0 & 0 & 0 & 0 \\
 0&0&0&0&0\\
  0 & a_2 & 0 & 0&0 \\
  a_1 & 0 & 0 & 0 &0\\
\end{array}%
\right). $$}
The coadjoint orbit  $\OSa$ is defined in \\
\\
{\small $\gx^* = \left\{\left(%
\begin{array}{ccccc}
  0 & 0 & 0 &0&0\\
  y_{21} & 0 & 0&0&0 \\
  y_{31} & y_{32} & 0&0&0\\
y_{41} & y_{42} & y_{43}&0&0\\
 y_{51} & y_{52} & y_{53}&y_{54}&0  \\
\end{array}%
\right) \right\}$} by the  equations ~~$ y_{51}=a_1,~~\left|
\begin{array}{cc}
  y_{41} & y_{42} \\
  y_{51} & y_{52} \\
\end{array}
\right| = -a_1a_2.$
\\
\\
 ii) ~$S = \{(1,4)\}$. Then  $L_S^+ =
\{(1,5)\}$, ~$L_S^-=\{(4,5)\}$,
 ~ $L_S^0=\{(2,5)\}$,~
$\La_S \cong  \Fbq^*\times \Fbq = \{(a_1,a_2)\in \Fbq^2:~ a_1\ne
0\}$,
{\small $$\px_S = \left\{\left(%
\begin{array}{ccccc}
  0 & x_{12} & x_{13} & x_{14}&x_{15}\\
  0 & 0 & x_{23}&0&x_{25} \\
  0 & 0 & 0 &0& 0\\
  0 & 0 & 0 &0& x_{45}\\
0 & 0 & 0 & 0&0\\
\end{array}%
\right)\right\}, \quad \la_{S,a} =
\left(%
\begin{array}{ccccc}
0&0&0&0&0\\
  0& 0 & 0 &0& 0 \\
  0 & 0 & 0 &0& 0 \\
  c_1 & 0 & 0 & 0&0 \\
  0 & a_1 &0& a_2 & 0 \\
\end{array}%
\right). $$} The coadjoint orbit  $\OSa$ is defined in $\gx^*$
by the equations $y_{51}=0$, ~$y_{52}=a_1$, ~$y_{41}=c$, ~ $ y_{52}y_{21} + y_{53}y_{31}+y_{54}y_{41} = a_2c_1 $. \\
iii) ~ $S = \{(2,3)\}$. Then $L_S^+ = \{(2,4)\}$,
~$L_S^-=\{(3,4)\}$,
 ~ $L_S^0=\{(1,5)\}$,~
$\La_S \cong  \Fbq^*\times \Fbq = \{(a_1,a_2)\in \Fbq^2:~ a_1\ne
0\}$,
{\small $$\px_S = \left\{\left(%
\begin{array}{ccccc}
  0 & x_{12} & x_{13} & x_{14}&x_{15}\\
  0 & 0 & x_{23}&x_{24}&0 \\
  0 & 0 & 0 &x_{34}& 0\\
  0 & 0 & 0 &0& 0\\
0 & 0 & 0 & 0&0\\
\end{array}%
\right)\right\}, \quad \la_{S,a} =
\left(%
\begin{array}{ccccc}
0&0&0&0&0\\
  0& 0 & 0 &0& 0 \\
  0 & c_2 & 0 &0& 0 \\
  0 & 0 & a_2 & 0&0 \\
  a_1 & 0 &0& 0 & 0 \\
\end{array}%
\right). $$} The coadjoint orbit  $\OSa$ is defined in $\gx^*$
by the  equations $$y_{51}=a_1,\quad \left|%
\begin{array}{cc}
  y_{41} & y_{42} \\
  y_{51} & y_{52} \\
\end{array}%
\right| = 0,\quad \left|%
\begin{array}{cc}
  y_{31} & y_{32} \\
  y_{51} & y_{52} \\
\end{array}%
\right| = -a_1c_2,\quad \left|%
\begin{array}{cc}
  y_{41} & y_{43} \\
  y_{51} & y_{53} \\
\end{array}%
\right| = -a_1a_2.$$ iv) ~ ~$S = \{(1,4), (2,3)\}$. Then $L_S^+ =
\{(1,5), (2,5)\}$, ~$L_S^-=\{(3,5), (4,5)\}$,
 ~ $L_S^0= \varnothing$,~
$\La_S \cong  \Fbq^2 = \{(a_1,a_2):~ a_1,a_2\in \Fbq\}$,
{\small $$\px_S = \left\{\left(%
\begin{array}{ccccc}
  0 & x_{12} & x_{13} & x_{14}&x_{15}\\
  0 & 0 & x_{23}&0&x_{25} \\
  0 & 0 & 0 &0& x_{35}\\
  0 & 0 & 0 &0& x_{45}\\
0 & 0 & 0 & 0&0\\
\end{array}%
\right)\right\}, \quad \la_{S,a} =
\left(%
\begin{array}{ccccc}
0&0&0&0&0\\
  0& 0 & 0 &0& 0 \\
  0 & c_2 & 0 &0& 0 \\
  c_1 & 0 & 0 & 0&0 \\
  0 & 0 &a_1& a_2 & 0 \\
\end{array}%
\right). $$} The coadjoint orbit  $\OSa$ is defined in $\gx^*$ by
the equations $y_{51}=y_{52}=0$, ~$y_{53}=a_1$, ~$y_{41}=c_1$, ~ $
y_{53}y_{31}+y_{54}y_{41} = a_2c_1 $,~~
$ \left|%
\begin{array}{cc}
  y_{31} & y_{32} \\
  y_{41} & y_{42} \\
\end{array}%
\right|=-c_1c_2$.\\
 The coadjoint orbit  $\OSa$ is defined in $\gx^*$ by the equations \\
$y_{31}=0$,~ $y_{21}=c$,~ $y_{32}=a$. Any representation of
Gelfand-Graev type $V(\la)$ decomposes into a sum of irreducible
representations that correspond to  mentioned above orbits. The
number of irreducible components equals to  $(2q-1)^2$. The matrices
$X_{S,b}$ from theorem  \ref{th5} has the  form {\small $$
\left(%
\begin{array}{ccccc}
  1 & 0 &0&0& b_1 \\
  0 & 1 &0& 0&b_2 \\
  0 & 0 &1&0& b_2\\
0 & 0 & 0&1 &b_1\\
0&0&0&0&1\\
\end{array}%
\right), \quad
\left(%
\begin{array}{ccccc}
  1 & 0 &0&0& b_1 \\
  0 & 1 &0& 0&a_1 \\
  0 & 0 &1&0& 0\\
0 & 0 & 0&1 &b_1\\
0&0&0&0&1\\
\end{array}%
\right), \quad
\left(%
\begin{array}{ccccc}
  1 & 0 &0&0& a_1 \\
  0 & 1 &0& a_1b_1&b_1 \\
  0 & 0 &1&0& b_1\\
0 & 0 & 0&1 &0\\
0&0&0&0&1\\
\end{array}%
\right), \quad
~\left(%
\begin{array}{ccccc}
  1 & 0 &0&0& a_1 \\
  0 & 1 &0& a_1a_2&a_2 \\
  0 & 0 &1&0& 0\\
0 & 0 & 0&1 &0\\
0&0&0&0&1\\
\end{array}%
\right),$$} where $ b_1,b_2\in\Fbq^*$, ~$ a_1,a_2\in \Fbq.$

\end{document}